\documentclass[12pt,fleqn]{article}
\usepackage{amssymb,amsthm,epsfig}
\newtheorem{theorem}{Theorem}[section]
\newtheorem{lemma}[theorem]{Lemma}
\newtheorem{corollary}[theorem]{Corollary}
\newtheorem{definition}[theorem]{Definition}

\newcommand{\old}[1]{{}}

\title{The extremal values of the Wiener index of a tree with given degree sequence}
\author{ 
Hua
 Wang\\
\tt{hua@math.ufl.edu} \\
Department of Mathematics,\\
University of Florida,
Gainesville, FL 32611\\
}

\begin{document}
\maketitle{}

\begin{abstract}
The Wiener index of a graph is the sum of the distances between all pairs of 
vertices, it has been one of the main descriptors that correlate a
chemical compound's molecular graph with experimentally gathered
data regarding the compound's characteristics. 
In \cite{wien}, the tree that minimizes the Wiener index among trees of given maximal degree is studied. 
We characterize
trees that achieve the maximum and minimum Wiener index, given the number of vertices and 
the degree sequence.
\end{abstract}


\newpage

\section{Terminology}

All graphs in this paper will be finite, simple and undirected.
A {\em tree} $T=(V,E)$ is a connected, acyclic graph. $V(T)$ denotes the vertex set of a tree $T$. 
We refer to vertices of degree 1 of $T$ as {\em leaves}.
 The unique path
connecting two vertices $v,u$ in $T$ will be denoted by
$P_{T}(v,u)$.
For a tree $T$ and two vertices $v$, $u$ of $T$, the
{\em distance} $d_{T}(v,u)$ between them is the number of edges
on the path $P_{T}(v,u)$.
 For a vertex $v$ of $T$, define the {\em distance of $v$} as
$ g_{T}(v)=\sum_{u\in V(T)} d_{T}(v,u).$  Then $ \sigma
(T)=\frac{1}{2}\sum_{v\in V(T)} g_{T}(v)$ denotes the {\em
Wiener index } of $T$.

For any vertex $v\in V(T)$, let $d(v)$ denote the {\em degree} of $v$, i.e. the number of edges incident to $v$.
The {\em degree sequence} of a tree is the sequence of the degrees (in descending order) of the non-leaf vertices.

We call a tree $(T, r)$ {\em  rooted at the vertex $r$} (or just $T$
if it is clear what the root is) by
specifying a vertex $r\in V(T)$. The {\em  height} of a vertex
$v$ of a
rooted tree $T$ with root $r$ is $h_{T}(v)=d_{T}(r,v)$.
 
For any two different vertices
$u, v$ in a rooted tree $(T,r)$, we say that $v$ is a {\em
successor} of $u$ and $u$ is an {\em ancestor} of $v$ if $P_{T}(r,u) \subset P_{T}(r,v)$.
Furthermore, if $u$ and $v$ are adjacent to each other and
$d_{T}(r,u)=d_{T}(r,v)-1$, we say that $u$ is the {\em  parent} of
$v$ and $v$ is a {\em  child } of $u$.
Two vertices $u,v$ are {\em siblings} of each other if they share the same parent.
A subtree of a tree will
often be described by its vertex set.
For a vertex $v$ in a rooted tree $(T, r)$, we use $T(v)$ to denote the subtree induced by 
$v$ and all its successors.

\section{Introduction}

To introduce our main results, we define the {\em greedy tree} (with a given degree sequence) as follows:

\begin{definition}
Suppose the degrees of the non-leaf vertices are given, the greedy tree is achieved
by the following 'greedy algorithm':

i) Label the vertex with the largest degree as $v$ (the root);

ii) Label the neighbors of $v$ as $v_1, v_2, \ldots$, assign the largest degrees available to them such that
$d(v_1) \geq d(v_2) \geq \ldots $;

iii) Label the neighbors of $v_1$ (except $v$) as $v_{11}, v_{12}, \ldots$ such that they take all the largest degrees available and that
$d(v_{11}) \geq d(v_{12}) \geq \ldots $, then do the same for $v_2$, $v_3$, $\ldots$;

iv) Repeat (iii) for all the newly labelled vertices, always start with the neighbors of 
the labelled vertex with largest degree whose neighbors are not labelled yet.
\end{definition}

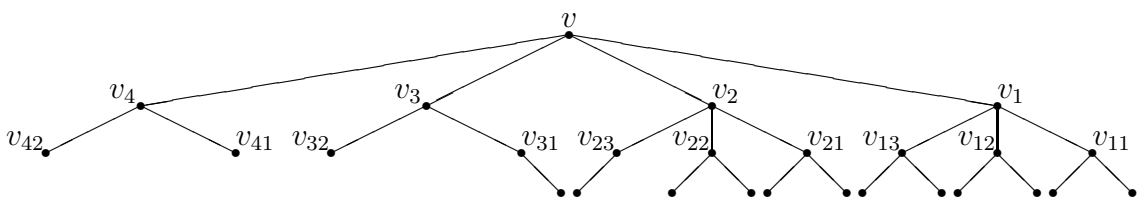
\begin{figure}[h]
\setlength{\unitlength}{3pt}
\begin{picture}(142,30)

\put(5,8){\circle*{1}}

\put(29,8){\circle*{1}}

\put(17,14){\circle*{1}} 

\put(17,14){\line(2,-1){12}} 
\put(17,14){\line(-2,-1){12}}

\put(41,8){\circle*{1}}

\put(65,8){\circle*{1}}

\put(53,14){\circle*{1}}

\put(65,8){\line(1,-1){5}}

\put(53,14){\line(2,-1){12}} 
\put(53,14){\line(-2,-1){12}}

\put(70,3){\circle*{1}} 

\put(72,3){\circle*{1}} 

\put(77,8){\circle*{1}}

\put(84,3){\circle*{1}} 
\put(94,3){\circle*{1}} 
\put(89,8){\circle*{1}}

\put(96,3){\circle*{1}} 
\put(106,3){\circle*{1}} 
\put(101,8){\circle*{1}}

\put(89,14){\circle*{1}}

\put(77,8){\line(-1,-1){5}}

\put(89,8){\line(1,-1){5}} 
\put(89,8){\line(-1,-1){5}}

\put(101,8){\line(1,-1){5}} 
\put(101,8){\line(-1,-1){5}}

\put(89,14){\line(2,-1){12}} 
\put(89,14){\line(-2,-1){12}}
\put(89,14){\line(0,-1){6}}

\put(108,3){\circle*{1}} 
\put(118,3){\circle*{1}} 
\put(113,8){\circle*{1}}

\put(120,3){\circle*{1}} 
\put(130,3){\circle*{1}} 
\put(125,8){\circle*{1}}

\put(132,3){\circle*{1}} 
\put(142,3){\circle*{1}} 
\put(137,8){\circle*{1}}

\put(125,14){\circle*{1}} 

\put(113,8){\line(1,-1){5}} 
\put(113,8){\line(-1,-1){5}}

\put(125,8){\line(1,-1){5}} 
\put(125,8){\line(-1,-1){5}}

\put(137,8){\line(1,-1){5}} 
\put(137,8){\line(-1,-1){5}}

\put(125,14){\line(2,-1){12}} 
\put(125,14){\line(-2,-1){12}}
\put(125,14){\line(0,-1){6}}

\put(71,23){\circle*{1}}
\put(71,23){\line(2,-1){18}}
\put(71,23){\line(-2,-1){18}}
\put(71,23){\line(6,-1){54}}
\put(71,23){\line(-6,-1){54}}

\put(70,24){\makebox{$v$}} 
\put(13,15){\makebox{$v_4$}}
\put(49,15){\makebox{$v_3$}}
\put(89,15){\makebox{$v_2$}}
\put(125,15){\makebox{$v_1$}}

\put(0,9){\makebox{$v_{42}$}}

\put(29,9){\makebox{$v_{41}$}}

\put(36,9){\makebox{$v_{32}$}}

\put(65,9){\makebox{$v_{31}$}}

\put(72,9){\makebox{$v_{23}$}}
\put(84,9){\makebox{$v_{22}$}}
\put(101,9){\makebox{$v_{21}$}}

\put(108,9){\makebox{$v_{13}$}}
\put(120,9){\makebox{$v_{12}$}}
\put(137,9){\makebox{$v_{11}$}}

\end{picture}

\caption{A greedy tree}
\label{greedy}
\end{figure}

Fig.~\ref{greedy} shows a greedy tree with degree sequence $\{ 4, 4, 4, 3,3,3,3,3,3,3,2,2 \}$.

From the definition of the greedy tree, we immediately get:

\begin{lemma}\label{gre}
A rooted tree $T$ with a given degree sequence is a greedy tree if:

i) the root $v$ has the largest degree;

ii) the heights of any two leaves differ by at most 1;

iii) for any two vertices $u$ and $w$, if $h_T(w) < h_T(v)$, then $d(w)\leq d(u)$;

iv) for any two vertices $u$ and $w$ of the same height, $d(u)> d(w) \Rightarrow d(u') \geq d(w')$ for any
successors $u'$ of $u$ and $w'$ of $w$ that are of the same height;

v) for any two vertices $u$ and $w$ of the same height, $d(u)> d(w) \Rightarrow d(u') \geq d(w') $ and 
$ d(u'') \geq d(w'') $ for any
siblings $u'$ of $u$ and $w'$ of $w$ or successors $u''$  of $u'$ and $w''$ of $w'$ of the same height.
\end{lemma}

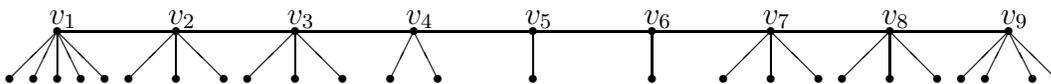
\begin{figure}[h]
\setlength{\unitlength}{3pt}
\begin{picture}(142,15)

\put(15,10){\line(1,0){120}} 

\put(15,10){\circle*{1}} 
\put(30,10){\circle*{1}} 
\put(45,10){\circle*{1}} 
\put(60,10){\circle*{1}} 
\put(75,10){\circle*{1}} 
\put(90,10){\circle*{1}} 
\put(105,10){\circle*{1}} 
\put(120,10){\circle*{1}} 
\put(135,10){\circle*{1}}  

\put(15,10){\line(1,-1){6}} 
\put(15,10){\line(-1,-1){6}}
\put(15,10){\line(1,-2){3}}
\put(15,10){\line(-1,-2){3}}
\put(15,10){\line(0,-1){6}}

\put(21,4){\circle*{1}} 
\put(18,4){\circle*{1}} 
\put(12,4){\circle*{1}} 
\put(9,4){\circle*{1}}  
\put(15,4){\circle*{1}}

\put(135,10){\line(1,-1){6}} 
\put(135,10){\line(-1,-1){6}}
\put(135,10){\line(1,-2){3}}
\put(135,10){\line(-1,-2){3}}

\put(141,4){\circle*{1}} 
\put(138,4){\circle*{1}} 
\put(132,4){\circle*{1}} 
\put(129,4){\circle*{1}}  

\put(30,10){\line(1,-1){6}} 
\put(30,10){\line(-1,-1){6}}
\put(30,10){\line(0,-1){6}}

\put(36,4){\circle*{1}} 
\put(24,4){\circle*{1}} 
\put(30,4){\circle*{1}}  

\put(120,10){\line(1,-1){6}} 
\put(120,10){\line(-1,-1){6}}
\put(120,10){\line(0,-1){6}}

\put(126,4){\circle*{1}} 
\put(114,4){\circle*{1}} 
\put(120,4){\circle*{1}}  

\put(45,10){\line(1,-1){6}} 
\put(45,10){\line(-1,-1){6}}
\put(45,10){\line(0,-1){6}}

\put(51,4){\circle*{1}} 
\put(39,4){\circle*{1}} 
\put(45,4){\circle*{1}}  

\put(105,10){\line(1,-1){6}} 
\put(105,10){\line(-1,-1){6}}
\put(105,10){\line(0,-1){6}}

\put(111,4){\circle*{1}} 
\put(99,4){\circle*{1}} 
\put(105,4){\circle*{1}}  

\put(60,10){\line(1,-2){3}} 
\put(60,10){\line(-1,-2){3}}

\put(63,4){\circle*{1}} 
\put(57,4){\circle*{1}} 

\put(75,10){\line(0,-1){6}}
\put(90,10){\line(0,-1){6}}

\put(75,4){\circle*{1}} 
\put(90,4){\circle*{1}}

\put(14,11){\makebox{$v_{1}$}}
\put(134,11){\makebox{$v_{9}$}}
\put(29,11){\makebox{$v_{2}$}}
\put(44,11){\makebox{$v_{3}$}}
\put(59,11){\makebox{$v_{4}$}}
\put(74,11){\makebox{$v_{5}$}}
\put(89,11){\makebox{$v_{6}$}}
\put(104,11){\makebox{$v_{7}$}}
\put(119,11){\makebox{$v_{8}$}}

\end{picture}

\caption{A greedy caterpillar}
\label{grepath}
\end{figure}

We also define the {\em greedy caterpillar} as a tree $T$ with given degree sequence \newline
$\{ d_1 \geq d_2 \geq \ldots \geq d_k \geq 2 \}$, 
that is formed by attaching pendant edges to a path $v_1v_2 \ldots v_k$ of length $k-1$ such that
$d(v_1) \geq d(v_k) \geq d(v_2) \geq d(v_{k-1}) \geq \ldots \geq d(v_{[\frac{k}{2} ]})$. Fig.~\ref{grepath} shows a 
greedy caterpillar with degree sequence $\{ 6, 5, 5, 5, 5, 5 ,4 ,3 , 3 \}$.

The structure of a chemical compound is usually modelled as a
polygonal shape, which is often called the \textit{molecular
graph} of this compound. The biochemical community has been using
topological indices to correlate a compound's molecular graph with
experimentally gathered data regarding the compound's
characteristics.

In 1947, Harold Wiener \cite{wiener} developed the Wiener Index. This concept has been one of the
most widely used descriptors in quantitative structure activity
relationships, as the Wiener index has been shown to have a strong
correlation with the chemical properties of the chemical compound.

Since the majority of the chemical applications of the Wiener
index deal with chemical compounds that have acyclic organic
molecules, whose molecular graphs are trees, the Wiener index of trees
has been extensively studied over the past years. 

It is well known that the
Wiener index is maximized by the path and minimized by the star among general trees
of the same size. Similar problems for more specific classes of trees seem to be more difficult.
In \cite{binl}, the Wiener index and the number of subtrees of binary trees
are studied, a not yet understood relation between them is discussed
for binary trees and trees in general. The correlation of various graph-theoretical indices including the Wiener index
is studied in the recent work of Wagner \cite{wagner}.

In \cite{wien}, the tree that minimizes the Wiener index among trees of given maximal degree is studied. 
However, the molecular graphs of the most practical interest have natural
restrictions on their degrees corresponding to the valences of the atoms, therefore 
it is reasonable to consider a tree with a fixed degree sequence. In this note, we study the extremal 
values of the Wiener index of a tree with given degree sequence and characterize these trees.
These trees are also shown to be the extremal trees with respect to {\em dominance order} by 
Fischermann, Rautenbach and Volkmann, for details see \cite{dom}.
We will prove the following:

\begin{theorem}\label{main}
Given the degree sequence and the number of vertices, the greedy tree minimizes the Wiener index.
\end{theorem}

\begin{theorem}\label{main'}
Given the degree sequence and the number of vertices, the greedy caterpillar maximizes the Wiener index.
\end{theorem}

In Section 3, a few Lemmas are given regarding the structure of an extremal tree with given degree sequence,
these results may be of interest on their own. 
We prove Theorem~\ref{main} in Section 4 and Theorem~\ref{main'} in Section 5.

\section{On the structure of an `optimal' tree}

For convenience, we will call a tree optimal if it minimizes the Wiener index among all trees with the same
number of vertices and the same degree sequence.

\begin{figure}[h]
\setlength{\unitlength}{3pt}
\begin{picture}(140,35)

\put(37,5){\line(1,0){61}} 
\put(75,5){\circle*{1}}
\put(90,5){\circle*{1}}
\put(60,5){\circle*{1}}
\put(45,5){\circle*{1}}
\put(67.5,5){\circle*{1}}

\put(31,5){\makebox{$\ldots$}} 
\put(100,5){\makebox{$\ldots$}} 

\put(74,7){\makebox{$x_1$}}
\put(89,7){\makebox{$x_2$}}
\put(59,7){\makebox{$y_1$}}
\put(44,7){\makebox{$y_2$}}
\put(66,7){\makebox{$z$}}

\put(74,15){\makebox{$X_1$}}
\put(89,15){\makebox{$X_2$}}
\put(59,15){\makebox{$Y_1$}}
\put(44,15){\makebox{$Y_2$}}
\put(66,15){\makebox{$Z$}}

\put(71,0){\dashbox{1}(10,20){}} 
\put(84,0){\dashbox{1}(12,20){}} 
\put(54,0){\dashbox{1}(10,20){}} 
\put(39,0){\dashbox{1}(12,20){}} 
\put(64,0){\dashbox{1}(7,20){}} 

\put(8,5){\line(1,0){21}} 
\put(21,5){\circle*{1}}
\put(10,5){\circle*{1}}

\put(106,5){\line(1,0){21}} 
\put(114,5){\circle*{1}}
\put(125,5){\circle*{1}}

\put(20,7){\makebox{$y_k$}}
\put(6,7){\makebox{$y_{k+1}$}}

\put(16,0){\dashbox{1}(11,20){}} 
\put(0,0){\dashbox{1}(14,20){}} 

\put(1,5){\makebox{$\ldots$}} 

\put(20,15){\makebox{$Y_k$}}
\put(2,15){\makebox{$Y_{>k}$}}

\put(113,7){\makebox{$x_k$}}
\put(123,7){\makebox{$x_{k+1}$}}

\put(108,0){\dashbox{1}(11,20){}} 
\put(122,0){\dashbox{1}(14,20){}} 

\put(129,5){\makebox{$\ldots$}} 

\put(113,15){\makebox{$X_k$}}
\put(127,15){\makebox{$X_{>k}$}}

\end{picture}

\caption{the components resulted from a path with $z$}
\label{path1}
\end{figure}
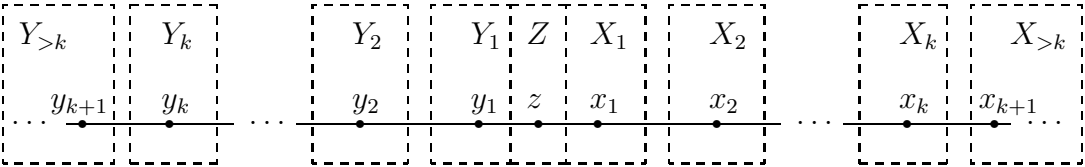

Consider a path in an optimal tree, after the removal of the edges on this path, some connected components will remain. Take 
a vertex and label it as $z$, and label the vertices on its right as $x_1, x_2, \ldots$, 
and the vertices on the left as $y_1, y_2, \ldots$. 
Let $X_i$ , $Y_i$ or $Z$ denote the component that contains the corresponding vertex. 
Let $X_{>k}$ and $Y_{>k}$ denote the trees induced by the vertices in 
$V(X_{k+1}) \cup V(X_{k+2}) \cup \ldots $ and $V(Y_{k+1}) \cup V(Y_{k+2}) \cup \ldots $ respectively (Fig.~\ref{path1}). 
Without loss of generality, assume
$|V(X_1)|\geq |V(Y_1)|$. 

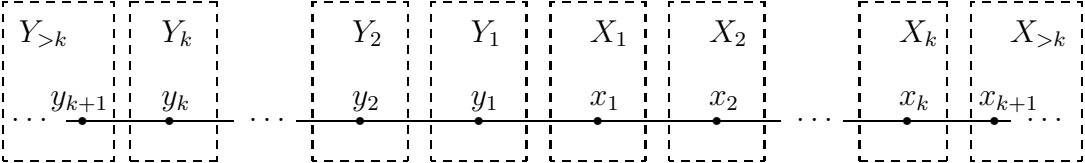
\begin{figure}[h]
\setlength{\unitlength}{3pt}
\begin{picture}(140,35)

\put(37,5){\line(1,0){61}} 
\put(75,5){\circle*{1}}
\put(90,5){\circle*{1}}
\put(60,5){\circle*{1}}
\put(45,5){\circle*{1}}

\put(31,5){\makebox{$\ldots$}} 
\put(100,5){\makebox{$\ldots$}} 

\put(74,7){\makebox{$x_1$}}
\put(89,7){\makebox{$x_2$}}
\put(59,7){\makebox{$y_1$}}
\put(44,7){\makebox{$y_2$}}

\put(74,15){\makebox{$X_1$}}
\put(89,15){\makebox{$X_2$}}
\put(59,15){\makebox{$Y_1$}}
\put(44,15){\makebox{$Y_2$}}

\put(69,0){\dashbox{1}(12,20){}} 
\put(84,0){\dashbox{1}(12,20){}} 
\put(54,0){\dashbox{1}(12,20){}} 
\put(39,0){\dashbox{1}(12,20){}} 

\put(8,5){\line(1,0){21}} 
\put(21,5){\circle*{1}}
\put(10,5){\circle*{1}}

\put(106,5){\line(1,0){21}} 
\put(114,5){\circle*{1}}
\put(125,5){\circle*{1}}

\put(20,7){\makebox{$y_k$}}
\put(6,7){\makebox{$y_{k+1}$}}

\put(16,0){\dashbox{1}(11,20){}} 
\put(0,0){\dashbox{1}(14,20){}} 

\put(1,5){\makebox{$\ldots$}} 

\put(20,15){\makebox{$Y_k$}}
\put(2,15){\makebox{$Y_{>k}$}}

\put(113,7){\makebox{$x_k$}}
\put(123,7){\makebox{$x_{k+1}$}}

\put(108,0){\dashbox{1}(11,20){}} 
\put(122,0){\dashbox{1}(14,20){}} 

\put(129,5){\makebox{$\ldots$}} 

\put(113,15){\makebox{$X_k$}}
\put(127,15){\makebox{$X_{>k}$}}

\end{picture}

\caption{the components resulted from a path without $z$}
\label{path}
\end{figure}

The next three lemmas hold for the path described above with (Fig.~\ref{path1}) or without (Fig.~\ref{path})
$z$.

\newpage

\begin{lemma}\label{geq1}
In the situation described above, if $|V(X_i)|\geq |V(Y_i)|$ for $i=1,2, \ldots, k$, then we can assume
\begin{equation}\label{geq}
|V(X_{>k})| \geq |V(Y_{>k})|
\end{equation}
in an optimal tree.
\end{lemma}
\begin{proof}
Suppose (for contradiction) that (\ref{geq}) does not hold. We will show that
switching $X_{>k}$ and $Y_{>k}$ (after which (\ref{geq}) holds) will not increase the Wiener index.

First, for a path without $z$,
note that in this operation, the lengths of the paths with both or neither end vertices in $V(X_{>k}) \cup V(Y_{>k})$
do not change. Hence we only need to consider the sum of the lengths of the paths that contain exactly one
end vertex in $V(X_{>k}) \cup V(Y_{>k})$.

For the distance between any vertex in $X_i$ ($i=1,2,\ldots, k$) and any vertex in $X_{>k}$,
this operation increases the distance by $2i-1$, then the total amount increased is 
\begin{equation}\label{incx}
\sum_{i=1}^k (2i-1) |V(X_i)| |V(X_{>k})|  ;
\end{equation}
Similarly, for the distances between any vertex in $Y_i$ ($i=1,2,\ldots, k$) 
and any vertex in $X_{>k}$, the total amount decreased is
\begin{equation}\label{decx}
\sum_{i=1}^k (2i-1) |V(Y_i)| |V(X_{>k})| ;
\end{equation}
For the distances between any vertex in $Y_i$ ($i=1,2,\ldots, k$) 
and any vertex in $Y_{>k}$, the total amount increased is
\begin{equation}\label{incy}
\sum_{i=1}^k (2i-1) |V(Y_i)| |V(Y_{>k})| ;
\end{equation}
For the distances between any vertex in $X_i$ ($i=1,2,\ldots, k$) 
and any vertex in $Y_{>k}$, the total amount decreased is
\begin{equation}\label{decy}
\sum_{i=1}^k (2i-1) |V(X_i)| |V(Y_{>k})| .
\end{equation}
Now $(\ref{incx}) + (\ref{incy}) - (\ref{decx}) - (\ref{decy})$
yields the total change of the Wiener index via this operation

$\sum_{i=1}^k (2i-1) (|V(X_i)| |V(X_{>k})| + |V(Y_i)| |V(Y_{>k})|
- |V(Y_i)| |V(X_{>k})| - |V(X_i)| |V(Y_{>k})| )$
\[
= \sum_{i=1}^k (2i-1) (|V(X_i)|-|V(Y_i)|)( |V(X_{>k})| - |V(Y_{>k})| )
\leq 0 .\]

For a path with $z$, note that the distance of a path with at least one end vertex in $Z$ 
does not change during this operation.
Similar to the first case, the total change of the Wiener index via this operation is
\[ \sum_{i=1}^k (2i) (|V(X_i)|-|V(Y_i)|)( |V(X_{>k})| - |V(Y_{>k})| )
\leq 0 .\]

\end{proof}

\begin{lemma}\label{geq2} 
If $|V(X_i)|\geq |V(Y_i)|$ for $i=1,2, \ldots, k-1$ and $|V(X_{>k})| \geq |V(Y_{>k})|$, then we can assume
\begin{equation}\label{ngeq}
|V(X_k)| \geq |V(Y_k)|
\end{equation}
in an optimal tree.
\end{lemma}
\begin{proof}
Suppose (for contradiction) that (\ref{ngeq}) does not hold, we will show that switching $X_{k}$ and $Y_k$ 
(after which (\ref{ngeq}) holds) will not 
increase the Wiener index.

Similar to the proof of Lemma~\ref{geq1}, the total change of the Wiener index via this operation is
\[ \sum_{i=1}^{k-1} (2i-1) (|V(X_i)|-|V(Y_i)|)( |V(X_{k})| - |V(Y_{k})| ) \]
\[ + 
(2k-1) (|V(X_{>k})|-|V(Y_{>k})|)( |V(X_{k})| - |V(Y_{k})| ) \leq 0 \]
for a path without $z$ and 
\[ \sum_{i=1}^{k-1} (2i) (|V(X_i)|-|V(Y_i)|)( |V(X_{k})| - |V(Y_{k})| ) \]
\[ + 
(2k) (|V(X_{>k})|-|V(Y_{>k})|)( |V(X_{k})| - |V(Y_{k})| ) \leq 0 \]
for a path with $z$.
\end{proof}

\begin{corollary}\label{geq3}
If $|V(X_i)|\geq |V(Y_i)|$ for $i=1,2, \ldots, k-1$ and $|V(X_{>k-1})| \geq |V(Y_{>k-1})|$, then we can assume
$d(x_k) \geq d(y_k)$ in an optimal tree.
\end{corollary}
\begin{proof}
Suppose (for contradiction) that $a=d(x_k) < d(y_k)=a+b$, 
the removal of $x_k$ ($y_k$) from $X_k$ ($Y_k$) will result in $a$ ($a+b$) components, take any $b$ components
(Let $B$ be the set of vertices in these $b$ components)
from $Y_k$ and attach them to $x_k$ (after which we have $d(x_k) \geq d(y_k)$) 
will preserve the degree sequence of the tree.

We will show that this operation will not increase the Wiener index.

Similar to the previous proofs, the total change of the Wiener index in this operation is
\[ \sum_{i=1}^{k-1} (2i-1) (|V(Y_i)|-|V(X_i)|)|B| \]
\[ + (2k-1) (|V(Y_{>k-1})|-|B|-|V(X_{>k-1})|)|B| \leq 0 \]
for a path without $z$ and 
\[ \sum_{i=1}^{k-1} (2i) (|V(Y_i)|-|V(X_i)|)|B| + 
(2k) (|V(Y_{>k-1})|-|B|-|V(X_{>k-1})|)|B| \leq 0 \]
for a path with $z$.
\end{proof}

\noindent {\bf Remark:}
In Lemmas~\ref{geq1}, \ref{geq2} and Corollary~\ref{geq3}, if at least one strict inequality holds in the conditions, then
the conclusion is forced and we can replace `can assume' by `must have' in the statement.

\newpage

Now, for a maximal path in an optimal tree, we can label the vertices and components with
vertices labelled as $w_1, w_2, \ldots $ and $u_1, u_2, \ldots $ and the components labelled as 
$W_i$ and $U_i$, while $U_1$ is the component with most vertices (Fig.~\ref{path2}) s.t.
the following hold:

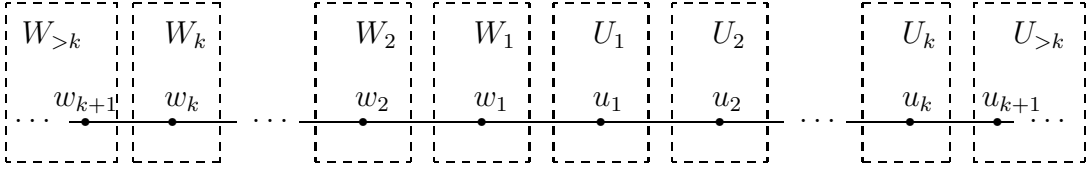
\begin{figure}[h]
\setlength{\unitlength}{3pt}
\begin{picture}(140,25)

\put(37,5){\line(1,0){61}} 
\put(75,5){\circle*{1}}
\put(90,5){\circle*{1}}
\put(60,5){\circle*{1}}
\put(45,5){\circle*{1}}

\put(31,5){\makebox{$\ldots$}} 
\put(100,5){\makebox{$\ldots$}} 

\put(74,7){\makebox{$u_1$}}
\put(89,7){\makebox{$u_2$}}
\put(59,7){\makebox{$w_1$}}
\put(44,7){\makebox{$w_2$}}

\put(74,15){\makebox{$U_1$}}
\put(89,15){\makebox{$U_2$}}
\put(59,15){\makebox{$W_1$}}
\put(44,15){\makebox{$W_2$}}

\put(69,0){\dashbox{1}(12,20){}} 
\put(84,0){\dashbox{1}(12,20){}} 
\put(54,0){\dashbox{1}(12,20){}} 
\put(39,0){\dashbox{1}(12,20){}} 

\put(8,5){\line(1,0){21}} 
\put(21,5){\circle*{1}}
\put(10,5){\circle*{1}}

\put(106,5){\line(1,0){21}} 
\put(114,5){\circle*{1}}
\put(125,5){\circle*{1}}

\put(20,7){\makebox{$w_k$}}
\put(6,7){\makebox{$w_{k+1}$}}

\put(16,0){\dashbox{1}(11,20){}} 
\put(0,0){\dashbox{1}(14,20){}} 

\put(1,5){\makebox{$\ldots$}} 

\put(20,15){\makebox{$W_k$}}
\put(2,15){\makebox{$W_{>k}$}}

\put(113,7){\makebox{$u_k$}}
\put(123,7){\makebox{$u_{k+1}$}}

\put(108,0){\dashbox{1}(11,20){}} 
\put(122,0){\dashbox{1}(14,20){}} 

\put(129,5){\makebox{$\ldots$}} 

\put(113,15){\makebox{$U_k$}}
\put(127,15){\makebox{$U_{>k}$}}

\end{picture}

\caption{the components resulted from a path}
\label{path2}
\end{figure}

\begin{lemma}\label{vcon}
In an optimal tree, we can 
label the vertices such that
\[ |V(U_1)| \geq |V(W_1)| \geq |V(U_2)| \geq |V(W_2)| \geq \ldots \geq |V(U_m)| = |V(W_m)| =1 \]
if the path is of odd  lenth ($2m-1$); and
\[ |V(U_1)| \geq |V(W_1)| \geq |V(U_2)| \geq |V(W_2)| \geq \ldots \geq |V(W_m)| = |V(U_{m+1})| =1 \]
if the path is of even  lenth ($2m$).
\end{lemma}
\begin{proof}
We only show the proof for a path of odd length, the other case is similar.
First,  we can assume $|V(U_1)| \geq |V(W_1)| \geq |V(U_2)|$ by symmetry. Now suppose we have
\begin{equation}\label{assum}
|V(U_1)| \geq |V(W_1)| \geq |V(U_2)| \geq |V(W_2)| \geq \ldots \geq |V(W_{k-1})| \geq |V(U_k)| 
\end{equation}
for some $k$. 

If equality holds in (\ref{assum}) except the last one, we can simply switch the label of $U_i$ and $W_i$ (if necessary)
to guarantee that $|V(U_k)| \geq |V(W_k)|$. Otherwise,
(\ref{assum}) implies $|V(U_{>k-1})| \geq |V(W_{>k-1})|$ by Lemma~\ref{geq1}, if 
$|V(W_k)| > |V(U_k)|$, then 
\[ |V(U_{>k})|=|V(U_{>k-1})| - |V(U_k)| > |V(W_{>k-1})|-|V(W_k)| = |V(W_{>k})|. \]
Applying Lemma~\ref{geq2} to $U_k$ and $W_k$ (in the setting that $x_i=u_i, y_i=w_i$ for $i=1,2,\ldots $)
yields a contradiction. Thus we have
\[ |V(U_1)| \geq |V(W_1)| \geq |V(U_2)| \geq |V(W_2)| \geq \ldots \geq |V(U_{k})| \geq |V(W_{k})|. \]
If all the equalities hold, we can switch the label of $U_{i+1}$ and $W_i$ for $i\geq 1$ (if necessary) and 
guarantee that $|V(W_k)|\geq |V(U_{k+1})|$. Otherwise,
apply Lemma~\ref{geq1} to $U_{>k}$ and $W_{>k-1}$ in the following setting:
\begin{equation}\label{set1} Z=U_1, Y_i=U_{i+1}, X_i=W_i, z=u_1, y_i=u_{i+1}, x_i=w_i,
\end{equation}
Then we have $|V(X_i)|\geq |V(Y_i)|$ for $i=1,2, \ldots , k-1$, 
thus 
\[ |V(W_{>k-1})|= |V(X_{>k-1})| \geq |V(Y_{>k-1})| = |V(U_{>k})| \] by Lemma~\ref{geq1}. If 
$|V(Y_k)|= |V(U_{k+1})| > |V(W_k)| = |V(X_k)|$, then 
\[ |V(Y_{>k})|=|V(Y_{>k-1})| - |V(Y_{k})| < |V(X_{>k-1})|-|V(X_k)| = |V(X_{>k})|. \]
Applying Lemma~\ref{geq2} to $Y_k=U_{k+1}$ and $X_k=W_k$ in setting (\ref{set1}) yields a contradiction. Thus we have
\[ |V(U_1)| \geq |V(W_1)| \geq |V(U_2)| \geq |V(W_2)| \geq \ldots \geq |V(U_{k})| \geq |V(W_{k})| \geq |V(U_{k+1})|. \]
The Lemma follows by induction.
\end{proof}

\noindent {\bf Remark:} 
Lemma~\ref{vcon} can be shown in a much easier way by using an equivalent definition of the Wiener index and 
simple application of a classic number theory result (\cite{hardy}). We keep the combinatorial proof here to 
provide a better understanding of the whole idea.

\begin{lemma}\label{dcon}
In an optimal tree, for a path with labelling as in Lemma~\ref{vcon}, we have
\[ d(u_1) \geq d(w_1) \geq d(u_2) \geq d(w_2) \geq \ldots \geq d(u_m) = d(w_m) = 1 \]
if the path is of odd  lenth ($2m-1$); and
\[ d(u_1) \geq d(w_1) \geq d(u_2) \geq d(w_2) \geq \ldots \geq d(u_m) \geq d(w_m) = d(w_{m+1}) =1 \]
if the path is of even  lenth ($2m$).
\end{lemma}
\begin{proof}
We only show the proof for the path of odd length, the other case is similar.

First, we have 
\[ |V(U_1)| \geq |V(W_1)| \geq |V(U_2)| \geq |V(W_2)| \geq \ldots \geq |V(U_m)| = |V(W_m)| = 1. \]
Now apply Corollary~\ref{geq3} to $u_i , u_{i+1}$ for $i=1,2, \ldots m-1$ in the following setting: 
\[ y_1=u_{i+1}, y_2=u_{i+2}, \ldots ; x_1=u_i, x_2=u_{i-1}, \ldots x_i=u_1, x_{i+1}=v_1, \ldots \]
Then $|V(X_{>1})|=\sum_{k=1}^{m} |V(W_k)| + \sum_{k=1}^{i-1} |V(U_k)| > \sum_{k=i+2}^{m} |V(U_k)| =|V(Y_{>1})|$,
implying that $d(u_i)=d(x_1) \geq d(y_1) = d(u_{i+1})$.

Thus we have
\[ d(u_1) \geq d(u_2) \geq \ldots \geq d(u_m). \]
Similarly, applying Corollary~\ref{geq3} to $w_i , w_{i+1}$ for $i=1,2, \ldots m-1$ yields
\[ d(w_1) \geq d(w_2) \geq \ldots \geq d(w_m). \]

For $u_i$ and $w_i$, if equality holds everywhere in Lemma~\ref{vcon}, we can 
again switch the labels and have $d(u_i) \geq d(w_i)$. 
Otherwise, applying Corollary~\ref{geq3} to $u_i , w_i$ (in the setting that $x_i=u_i, y_i=w_i$ for $i=1,2, \ldots $) 
yields that $d(u_i) \geq d(w_i)$ for $i=1,2, \ldots , m$;

Similarly, applying Corollary~\ref{geq3} to $w_i, u_{i+1}$ in the setting (\ref{set1}) 
yields that $d(w_i) \geq d(u_{i+1})$ for $i=1,2, \ldots, m-1$.
\end{proof}

\section{Proof of Theorem~\ref{main}}

It has been shown that $g_T(v)$ is minimized at one or two adjacent vertices on any path and hence in the 
whole tree (called the {\em centroid} of the tree), 
see \cite{jordan} and \cite{zelinka} for details. From Lemma~\ref{vcon} and Lemma~\ref{dcon}, simple calculation shows:

On any path of an optimal tree labelled as in Lemma~\ref{vcon} and Lemma~\ref{dcon}, 
\begin{equation}\label{cor}
\hbox{ the minimal value of $g_T(v)$ is achieved at
$u_1$ }
\end{equation}
where $d(u_1)$ and $|V(U_1)|$ are maximum on the path.

There are two cases:

\noindent i) If there is only one vertex in the centroid, label it as $v$.

\noindent ii) If there are two vertices in the centroid, label the vertex in the component consisting of more vertices 
(after the removal of the edge in between the two vertices in the centroid) as $v$ and the other one as $v_1$. If the two 
components contain the same number of vertices, just choose either one as $v$ and the other one as $v_1$.

\begin{proof} 
We only show the first case, the second one is similar.

In an optimal tree $T$, consider $T$ as rooted at $v$, we know $v$ is of the largest degree 
immediately from (\ref{cor}) (hence (i) of Lemma~\ref{gre} is satisfied). 

Consider any path starting at a leaf $u$, passing $v$, ending at a leaf $w$ whose only common ancestor with $u$ is $v$.
Apply Lemma~\ref{dcon} to this path such that $u_1=v$, we must have
$|d_T(u, v) - d_T(w, v)|\leq 1$, then the heights of any two leaves differ by at most 1 
(hence (ii) of Lemma~\ref{gre} is satisfied). 
Furthermore, it is also implied that 
\begin{equation}\label{neweq}
\hbox{$d(x)\geq d(y)$ for any two
vertices such that $y$ is a successor of $x$.}
\end{equation}
For vertex $x$ of height $i$ and vertex $y$ of height $j$ ($i< j$), consider the following two cases:

a) if $y$ is a successor of $x$, then we have $d(x)\geq d(y)$ from (\ref{neweq});

b) otherwise, let $u$ be the common ancestor of them
that is on the path $P_T(x,y)$, apply Lemma~\ref{dcon} to the path that passes through $y', y, u, x, x'$,
where $y',x'$ are leaves that are successors (or equal to) $y,x$ respectively. 
We must have $u_1=u$ by (\ref{cor}) and Lemma~\ref{vcon}, then $x=u_{k+1}, y=w_{l}$ or $x=w_{k}, y=u_{l+1}$, where 
$k=i-h_T(u), l=j-h_T(u), k+1\leq l$. Either way, Lemma~\ref{dcon} implies that
$d(x)\geq d(y)$.

Hence (iii) of Lemma~\ref{gre} is satisfied.

For two non-leaf vertices $x$ and $y$ of the same height $i$ with $d(x) > d(y)$, let $x'$ and
$y'$ (of the same height $j$) be the successors of $x$ and $y$ respectively. Apply 
Lemma~\ref{dcon} to the longest path that passes through $y', y, u, x, x'$, where
$u$ is the
common ancestor of $x,y$ that is on the path $P_T(x,y)$.
We must have $u_1=u$ by (\ref{cor}) and Lemma~\ref{vcon}, 
then $x=w_{k}, x'=w_{l}, y=u_{k+1}, y'=u_{l+1}$ as $d(x) > d(y)$,
where $k= i-h_T(u), l=j-h_T(u)$.
Thus implying $d(x') \geq d(y')$ 
(hence (iv) of Lemma~\ref{gre} is satisfied). 

Now let $x_0$ ($x'$) and $y_0$ ($y'$) be the parents (siblings) of $x$ and $y$ respectively, 
let $x''$ and $y''$ (of the same height $j$)
be successors of $x'$ and $y'$ respectively.  The conclusion of (iv) implies
\begin{equation}\label{local} 
|V(T(x_0)/T(x'))| > |V(T(y_0)/T(y'))|.
\end{equation}
Now consider the longest path that passes through $y'', y', u, x', x''$, where
$u$ is the common ancestor of $x$ and $y$ that is on the path $P_T(x', y')$.
Apply Lemma~\ref{dcon},
we must have $u_1=u$ by (\ref{cor}) and Lemma~\ref{vcon}, then $x'=w_{k}, x''=w_{l}, y'=u_{k+1}, y''=u_{l+1}$
by (\ref{local}) and Lemma~\ref{vcon}, 
where $k= i-h_T(u), l=j-h_T(u)$.
Thus we have  $d(x') \geq d(y')$ and $d(x'') \geq d(y'')$ (hence (v) of Lemma~\ref{gre} is satisfied).

In conclusion, by Lemma~\ref{gre}, the optimal tree is the greedy tree.
\end{proof}

\section{On Theorem~\ref{main'} }

Similar to Lemma~\ref{vcon}, we have the following for trees with given degree sequence that
maximize the Wiener index (refer to Fig.~\ref{path2}), we leave the proof to the reader:

\begin{lemma}\label{vcon'}
In a tree 
with a given number of vertices and degree sequence that maximizes the Wiener index, we can label
the vertices on the path with $U_1$ being the component consisting of the least vertices such that:
\[ |V(U_1)| \leq |V(W_1)| \leq |V(U_2)| \leq |V(W_2)| \leq \ldots \leq |V(U_{m-1})| \leq |V(W_{m-1})| \]
if the path is of odd  lenth ($2m-1$); and
\[ |V(U_1)| \leq |V(W_1)| \leq |V(U_2)| \leq |V(W_2)| \leq \ldots \leq |V(W_{m-1})| \leq |V(U_{m})| \]
if the path is of even  lenth ($2m$).
\end{lemma}

\begin{proof}(of Theorem~\ref{main'})

Let $T$ be the tree that maximizes the Wiener index with a given degree sequence. Consider the longest path,
without loss of generality, let the path be $w_m w_{m-1} \ldots w_1 u_1 u_2 \ldots u_m$ of odd length (the other case
is similar).

First we show that every vertex not on the path is a leaf, otherwise, 
let $x$ be a neighbor of $w_i$ (the case for $u_i$ is similar) that is not on the path and is not a leaf.  
Consider the longest path that contains $w_m, w_i, x$, i.e. 
$w_m \ldots w_i x x_1 \ldots x_s y$ where $y$ is a leaf.

Let $W_i, U_i$ denote the components with respect to the path $w_m w_{m-1} \ldots w_1 u_1 u_2 \ldots u_m$
as in Lemma~\ref{vcon'}. Let $X_{w_m}, X_{w_{m-1}}, \ldots, X_{w_i}, X_{x}, X_{x_1}, \ldots, X_{x_s}$ denote the components
resulting from removing the edges on the path $w_m \ldots w_i x x_1 \ldots x_s y$. 
Now consider the path $w_m \ldots w_i x x_1 \ldots x_s y$, we have
\[ |V(X_{w_i})| \geq |V(U_{m-1})| \geq |V(W_i)| > |V(X_{x_s})|, \]
contradicting to Lemma~\ref{vcon'} (note that $i\leq m-2$). 
Thus, for every vertex on the path $w_m w_{m-1} \ldots w_1 u_1 u_2 \ldots u_m$, 
if it has any neighbor that is not on the path, they must be
leaves. Applying Lemma~\ref{vcon'} to the path $w_m w_{m-1} \ldots w_1 u_1 u_2 \ldots u_m$
yields that $T$ must be a greedy caterpillar.
\end{proof}

\end{document}